\documentclass[12pt]{amsart} 
\usepackage{amssymb,amscd}
\usepackage{verbatim}

\begin{document}

\numberwithin{equation}{section}

\newtheorem{theorem}[equation]{Theorem}
\newtheorem{lemma}[equation]{Lemma}
\newtheorem{conjecture}[equation]{Conjecture}
\newtheorem{proposition}[equation]{Proposition}
\newtheorem{corollary}[equation]{Corollary}
\newtheorem{cor}[equation]{Corollary}

\theoremstyle{definition}
\newtheorem*{definition}{Definition}
\newtheorem{example}[equation]{Example}

\theoremstyle{remark}
\newtheorem{remark}[equation]{Remark}
\newtheorem{remarks}[equation]{Remarks}
\newtheorem*{acknowledgement}{Acknowledgements}


\newenvironment{notation}[0]{%
  \begin{list}%
    {}%
    {\setlength{\itemindent}{0pt}
     \setlength{\labelwidth}{4\parindent}
     \setlength{\labelsep}{\parindent}
     \setlength{\leftmargin}{5\parindent}
     \setlength{\itemsep}{0pt}
     }%
   }%
  {\end{list}}

\newenvironment{parts}[0]{%
  \begin{list}{}%
    {\setlength{\itemindent}{0pt}
     \setlength{\labelwidth}{1.5\parindent}
     \setlength{\labelsep}{.5\parindent}
     \setlength{\leftmargin}{2\parindent}
     \setlength{\itemsep}{0pt}
     }%
   }%
  {\end{list}}
\newcommand{\Part}[1]{\item[\upshape#1]}

\renewcommand{\a}{\alpha}
\renewcommand{\b}{\beta}
\newcommand{\g}{\gamma}
\renewcommand{\d}{\delta}
\newcommand{\e}{\epsilon}
\newcommand{\f}{\phi}
\renewcommand{\l}{\lambda}
\renewcommand{\k}{\kappa}
\newcommand{\lhat}{\hat\lambda}
\newcommand{\m}{\mu}
\renewcommand{\o}{\omega}
\renewcommand{\r}{\rho}
\newcommand{\rbar}{{\bar\rho}}
\newcommand{\s}{\sigma}
\newcommand{\sbar}{{\bar\sigma}}
\renewcommand{\t}{\tau}
\newcommand{\z}{\zeta}

\newcommand{\D}{\Delta}

\newcommand{\gp}{{\mathfrak{p}}}
\newcommand{\gP}{{\mathfrak{P}}}
\newcommand{\gq}{{\mathfrak{q}}}
\newcommand{\gf}{{\mathfrak{f}}}

\newcommand{\Acal}{{\mathcal A}}
\newcommand{\Bcal}{{\mathcal B}}
\newcommand{\Ccal}{{\mathcal C}}
\newcommand{\Dcal}{{\mathcal D}}
\newcommand{\Ecal}{{\mathcal E}}
\newcommand{\F}{{\mathcal F}}
\newcommand{\Fcal}{{\mathcal F}}
\newcommand{\cF}{{\mathcal F}}
\newcommand{\Gcal}{{\mathcal G}}
\newcommand{\Hcal}{{\mathcal H}}
\newcommand{\Ical}{{\mathcal I}}
\newcommand{\Jcal}{{\mathcal J}}
\newcommand{\Kcal}{{\mathcal K}}
\newcommand{\Lcal}{{\mathcal L}}
\newcommand{\Mcal}{{\mathcal M}}
\newcommand{\Ncal}{{\mathcal N}}
\newcommand{\Ocal}{{\mathcal O}}
\newcommand{\Pcal}{{\mathcal P}}
\newcommand{\Qcal}{{\mathcal Q}}
\newcommand{\Rcal}{{\mathcal R}}
\newcommand{\Scal}{{\mathcal S}}
\newcommand{\Tcal}{{\mathcal T}}
\newcommand{\Ucal}{{\mathcal U}}
\newcommand{\Vcal}{{\mathcal V}}
\newcommand{\Wcal}{{\mathcal W}}
\newcommand{\Xcal}{{\mathcal X}}
\newcommand{\Ycal}{{\mathcal Y}}
\newcommand{\Zcal}{{\mathcal Z}}

\newcommand{\OO}{{\mathcal O}}    
\newcommand{\KK}{{\mathcal K}}    
\renewcommand{\O}{{\mathcal O}}   

\renewcommand{\AA}{\mathbb{A}}
\newcommand{\BB}{\mathbb{B}}
\newcommand{\CC}{\mathbb{C}}
\newcommand{\FF}{\mathbb{F}}
\newcommand{\GG}{\mathbb{G}}
\newcommand{\PP}{\mathbb{P}}
\newcommand{\NN}{\mathbb{N}}
\newcommand{\QQ}{\mathbb{Q}}
\newcommand{\RR}{\mathbb{R}}
\newcommand{\ZZ}{\mathbb{Z}}

\newcommand{\bfa}{{\mathbf a}}
\newcommand{\bfb}{{\mathbf b}}
\newcommand{\bfc}{{\mathbf c}}
\newcommand{\bfe}{{\mathbf e}}
\newcommand{\bff}{{\mathbf f}}
\newcommand{\bfg}{{\mathbf g}}
\newcommand{\bfp}{{\mathbf p}}
\newcommand{\bfr}{{\mathbf r}}
\newcommand{\bfs}{{\mathbf s}}
\newcommand{\bft}{{\mathbf t}}
\newcommand{\bfu}{{\mathbf u}}
\newcommand{\bfv}{{\mathbf v}}
\newcommand{\bfw}{{\mathbf w}}
\newcommand{\bfx}{{\mathbf x}}
\newcommand{\bfy}{{\mathbf y}}
\newcommand{\bfz}{{\mathbf z}}
\newcommand{\bfA}{{\mathbf A}}
\newcommand{\bfB}{{\mathbf B}}
\newcommand{\bfC}{{\mathbf C}}
\newcommand{\bfF}{{\mathbf F}}
\newcommand{\bfG}{{\mathbf G}}
\newcommand{\bfI}{{\mathbf I}}
\newcommand{\bfM}{{\mathbf M}}
\newcommand{\bfzero}{{\boldsymbol{0}}}
\newcommand{\bfmu}{{\boldsymbol\mu}}

\def\ta{{\tilde{a}}}
\def\tA{{\tilde{A}}}
\def\tb{{\tilde{b}}}
\def\tc{{\tilde{c}}}
\def\td{{\tilde{d}}}
\def\tf{{\tilde{f}}}
\def\tF{{\tilde{F}}}
\def\tg{{\tilde{g}}}
\def\tG{{\tilde{G}}}
\def\th{{\tilde{h}}}
\def\tK{{\tilde{K}}}
\def\tk{{\tilde{k}}}
\def\tz{{\tilde{z}}}
\def\tw{{\tilde{w}}}
\def\tphi{{\tilde{\varphi}}}
\def\tdelta{{\tilde{\delta}}}
\def\tbeta{{\tilde{\beta}}}
\def\talpha{{\tilde{\alpha}}}
\def\ttheta{{\tilde{\theta}}}
\def\tB{{\tilde{B}}}
\def\tR{{\tilde{R}}}
\def\tT{{\tilde{T}}}
\def\tE{{\tilde{E}}}
\def\tU{{\tilde{U}}}
\def\tGamma{{\tilde{\Gamma}}}

\newcommand{\ab}{{\textup{ab}}}
\newcommand{\Ahat}{{\hat A}}
\newcommand{\Aut}{\operatorname{Aut}}
\newcommand{\Cond}{{\mathfrak{N}}} 
\newcommand{\Disc}{\operatorname{Disc}}
\newcommand{\Div}{\operatorname{Div}}
\newcommand{\End}{\operatorname{End}}
\newcommand{\Frob}{\operatorname{Frob}}
\newcommand{\Fpbar}{{\overline{\FF}_p}}
\newcommand{\GK}{G_{\Kbar/K}}
\newcommand{\GL}{\operatorname{GL}}
\newcommand{\Gal}{\operatorname{Gal}}
\newcommand{\hhat}{{\hat h}}
\newcommand{\Image}{\operatorname{Image}}
\newcommand{\into}{\hookrightarrow}     
\newcommand{\Kbar}{{\bar K}}
\newcommand{\Kvbar}{{{\bar K}_v}}
\newcommand{\MOD}[1]{~(\textup{mod}~#1)}
\newcommand{\Norm}{\operatorname{N}}
\newcommand{\notdivide}{\nmid}
\newcommand{\nr}{{\textup{nr}}}    
\newcommand{\ord}{\operatorname{ord}}
\newcommand{\Pic}{\operatorname{Pic}}
\newcommand{\Qbar}{{\bar{\QQ}}}
\newcommand{\Qpbar}{{\bar{\QQ}_p}}
\newcommand{\kbar}{{\bar{k}}}
\newcommand{\QQbar}{{\bar{\QQ}}}
\newcommand{\rank}{\operatorname{rank}}
\newcommand{\res}{\operatornamewithlimits{res}}
\newcommand{\Resultant}{\operatorname{Res}}
\newcommand{\Res}{\operatorname{Res}}
\renewcommand{\setminus}{\smallsetminus}
\newcommand{\Spec}{\operatorname{Spec}}
\newcommand{\PGL}{\operatorname{PGL}}
\newcommand{\Support}{\operatorname{Supp}}
\newcommand{\tors}{{\textup{tors}}}
\newcommand{\val}{{\operatorname{val}}}
\newcommand{\<}{\langle}
\newcommand{\la}{{\langle}}
\renewcommand{\>}{\rangle}
\newcommand{\ra}{{\rangle}}
\newcommand{\Berk}{{\rm Berk}}
\newcommand{\Rat}{{\rm Rat}}
\newcommand{\RL}{{\rm RL}}
\newcommand{\HH}{{\mathbf H}}
\newcommand{\Gm}{{\mathbf G}_m}
\newcommand{\Hhat}{{\hat{H}}}


\hyphenation{archi-me-dean}


\title[Average values of dynamical Green's functions]{A lower bound for average values of dynamical Green's functions}

\author[Matt Baker]{Matthew Baker}
\email{mbaker@math.gatech.edu}
\address{Department of Mathematics,
          Georgia Institute of Technology, Atlanta GA 30332-0160, USA}

\begin{abstract}
We provide a Mahler/Elkies-style lower bound for the average values of dynamical Green's functions on $\PP^1$ over an arbitrary
valued field, and give some dynamical and arithmetic applications.
\end{abstract}

\thanks{The author would like to thank Rob Benedetto for helpful conversations while this work was in progress.  He would also
like to thank Laura DeMarco, Xander Faber, Robert Rumely, and Joe Silverman for their feedback on an earlier version of this paper.  
The author's research was supported by NSF Research Grant DMS-0300784.}

\maketitle


\section{Introduction}

\subsection{Background}

Throughout this paper, $K$ will denote a field endowed with an absolute value, which can be either 
archimedean or non-archimedean.
Let $\varphi \in K(T)$ be a rational function of degree $d \geq 2$, and let
$g_{\varphi}$ be the normalized Arakelov-Green's function associated to $\varphi$ 
(see (\ref{gFvDef}) below for a definition).
In several different situations, one is led to consider lower bounds for $g_{\varphi}$-discriminant sums of the form
\[
D_\varphi(z_1,\ldots,z_N) = \sum_{\substack{1 \leq i,j \leq N \\ i \neq j}} g_{\varphi}(z_i,z_j).
\]
Since $g_{\varphi}(z,w)$ is bounded below, one knows that 
$D_\varphi(z_1,\ldots,z_N) \geq -C N^2$ for some constant $C>0$.
The main result of this paper is the following stronger estimate:

\begin{theorem}
\label{GreenTheorem}
There is an effective constant $C> 0$, depending on $\varphi$ and $K$, such that if 
$N \geq 2$ and $z_1,\ldots,z_N$ are distinct points of $\PP^1(K)$, then
\begin{equation}
\label{mainboundequation}
D_\varphi(z_1,\ldots,z_N) \geq -C N \log N.
\end{equation}
\end{theorem}

We will first give some historical context for this result, and then we will explain some of its
applications.

\medskip

The first estimate of this kind appeared in the work of K.~Mahler \cite{Mahler},
who proved a useful inequality between the discriminant and Mahler measure of a polynomial.  
Specifically, suppose that $f(x)=a_0x^N+a_1x^{N-1}+\cdots+a_N=a_0(x-\alpha_1)\cdots(x-\alpha_N)$ is a polynomial with complex coefficients and degree $N\geq 2$, 
let 
$\Disc(f)=a_0^{2N-2}\prod_{i < j}(\alpha_i-\alpha_j)^2$ 
be the discriminant of $f$, and let
$M(f)= |a_0|\prod_{j=1}^N \max(1,|\alpha_j|)$ be the Mahler measure of $f$.
Mahler's inequality states that 
\begin{equation}
\label{MahlerIneq}
|\Disc(f)|\leq N^N\{M(f)\}^{2N-2}.
\end{equation}
The inequality (\ref{MahlerIneq}) can be used, for example, to give
a short proof of Bilu's equidistribution theorem (see e.g. \cite{Bombieri}).
Mahler's proof of (\ref{MahlerIneq}) uses Hadamard's inequality to estimate the determinant of a suitable Vandermonde matrix.
Essentially the same argument gives rise to the following reformulation of (\ref{MahlerIneq}).
For $z,w \in \PP^1(\CC)$ with $z \neq w$, define 
\[
g(z,w) = \left\{ 
\begin{array}{ll} 
-\log |z-w| + \log^+|z| + \log^+|w|  & z,w \neq \infty \\
\log^+|z|  & w = \infty \\
\log^+|w|  & z = \infty. \\
\end{array}
\right.
\]
Then if $N \geq 2$ and $z_1,\ldots,z_N \in \PP^1(\CC)$ are distinct points, we have
\begin{equation}
\label{AltMahlerIneq}
\sum_{i \neq j} g(z_i,z_j) \geq -N\log N.
\end{equation}
Note that $g(z,w) \geq -\log 2$ for all $z\neq w$, and thus (\ref{AltMahlerIneq}) gives a significant improvement over the
trivial lower bound of $-(\log2) N(N-1)$.

\medskip

The reformulation (\ref{AltMahlerIneq}) of inequality (\ref{MahlerIneq}) is strongly reminiscent of the following result of Elkies 
(see \cite[\S VI, Theorem 5.1]{LangAT} and \cite[Appendix]{BP}).
Let $E/\CC$ be an elliptic curve,
and let $g(z,w)$ be the normalized Arakelov-Green's function associated to
the Haar measure $\mu$ on $E(\CC)$.
Then there exists a constant $C>0$ such that if $N \geq 2$ and $z_1,\ldots,z_N$ are distinct points of $E(\CC)$, 
then 
\begin{equation}
\label{ElkiesBound}
\sum_{i \neq j} g(z_i,z_j) \geq -C N \log N.
\end{equation}

The result proved in \cite{LangAT} is actually a generalization of (\ref{ElkiesBound}) to 
an arbitrary compact Riemann surface $X$ of genus at least 1, where Haar measure is replaced by
the canonical volume form on $X$.  
Lang remarks on page 152 that Elkies' argument can be used to prove a similar theorem for an arbitrary compact
Riemannian manifold.  

\medskip

Bounds of the form (\ref{ElkiesBound}) are
used in Arakelov theory to prove the existence of small sections of large powers of arithmetically ample metrized line bundles.
In addition, the bound (\ref{ElkiesBound}) is used in a series of papers by Hindry and Silverman
to obtain lower bounds for the canonical height of a non-torsion 
$k$-rational point on an elliptic curve $E$ over a number field $k$, as well as
upper bounds for the number of $k$-rational torsion points on $E$.
Further applications of (\ref{ElkiesBound}) to canonical heights on elliptic curves can be found in \cite{BP}.
The proof of (\ref{ElkiesBound}) for a general Riemann surface $X$
uses the spectral theory of the Laplacian operator; when $X=E$ is an elliptic curve, this amounts to Fourier
analysis on the complex torus $\CC / \Lambda$.
An analogue of (\ref{ElkiesBound}) for metrized graphs is proved in 
\cite{HAOMG}.

\medskip

Theorem~\ref{GreenTheorem} is a simultaneous generalization of (\ref{AltMahlerIneq}) and 
(\ref{ElkiesBound}) which applies to 
Arakelov-Green's functions attached to an arbitrary dynamical system on $\PP^1$.  It
works equally well over archimedean and non-archimedean fields.
The case where $\varphi(T)=T^2$ and $K=\CC$ corresponds to Mahler's result, and the case where $\varphi(T)$ is the degree 4 ``Latt{\`e}s map''
corresponding to multiplication by 2 on the quotient of an elliptic curve over $\CC$
by the involution $P \mapsto -P$ corresponds to (\ref{ElkiesBound}).
Since the analytic tools employed by Elkies do not seem to be available in the context of the dynamical systems attached to
arbitrary rational functions, our proof of (\ref{mainboundequation}) is based on a suitable modification of Mahler's original approach.  
In particular, we obtain a new and much more elementary proof of (\ref{ElkiesBound}), albeit with different constants.

\medskip

The possibility of proving (\ref{mainboundequation}) by a Mahler-style argument
was suggested to us by Lemma 4.1 of R.~Benedetto's paper \cite{Benedetto},
which gives a bound of the form (\ref{mainboundequation}) in the special case where 
$\varphi$ is a polynomial and $z_1,\ldots,z_N$ belong to the filled Julia set of $\varphi$.
Note that the case where $\varphi$ is a polynomial is simpler than the general case, due to the fact that 
polynomial maps have a superattracting fixed point at infinity.


\medskip

When $K$ is non-archimedean, the $-CN\log N$ term in (\ref{ElkiesBound}) can be replaced by $-CN$,
see \cite[Lemma 2.1]{HindrySilvermanIII}.
For general rational maps over a non-archimedean field, however, the bound (\ref{mainboundequation}) is optimal.
This can be seen using Lemma 3.4(a) and Remark 4.3 of \cite{Benedetto}.

\subsection{Definition and properties of $g_{\varphi}$}

In order to define the Arakelov-Green's function $g_{\varphi}$ which appears in the statement of
(\ref{mainboundequation}), 
we begin by recalling some definitions and results from \cite{BR}.

Write $\varphi$ in the form
\[
\varphi([z_0:z_1]) \ = \ [F_1(z_0,z_1):F_2(z_0,z_1)]
\]
for some homogeneous polynomials $F_1,F_2 \in K[x,y]$ of degree $d \geq 2$
with no common linear factor over $\Kbar$.
The polynomials $F_1,F_2$ are uniquely determined by $\varphi$ up to
multiplication by a common scalar $c \in K^*$.
The mapping
\[
F = (F_1,F_2) : \AA^2(K) \to \AA^2(K)
\]
is a lifting of $\varphi$ to $\AA^2$.

Let $\Res(F) := \Res(F_1,F_2)$ denote the homogeneous resultant of the
polynomials $F_1$ and $F_2$ (see \cite[\S 5.8]{VdW},\cite[\S 6]{DeMarco}).  
Since $F_1$ and $F_2$ have no common linear factor over $\Kbar$, 
we have $\Res(F) \neq 0$.

Let $\|(z_0,z_1)\| := \max \{ |z_0|,|z_1| \} $.
For $z = (z_0,z_1) \in K^2 \backslash \{ 0 \}$,
define the {\em homogeneous local dynamical height} $\Hhat_{F} : K^2
\backslash \{ 0 \} \to \RR$ by
\[
\Hhat_{F}(z) := \lim_{n\to\infty} \frac{1}{d^n} \log \| F^{(n)}(z) \|.
\]

By convention, we define $\Hhat_{F}(0,0) := -\infty$.
It is proved in \cite{BR} that the limit $\lim_{n\to\infty} \frac{1}{d^n} \log \| F^{(n)}(z) \|$
exists for all $z \in K^2 \backslash \{ 0 \}$, and that
$\frac{1}{d^n} \log \| F^{(n)}(z) \|$ converges uniformly on
$K^2 \backslash \{ 0 \}$ to $\Hhat_{F}(z)$.  (Note that the results in \cite{BR} are stated in the special case $K=\CC_v$, but
the proofs extend without modification to the more general case considered here.)
The definition of $\Hhat_{F}$ is independent of the
norm used to define it.  

\medskip

For $z = (z_0,z_1), w = (w_0, w_1) \in K^2$, put
\[
z \wedge w \ := \ z_0 w_1 - z_1 w_0 \ .
\]

When $z, w \in K^2$ are linearly independent over $K$, we define
\[
G_{F}(z,w) \ := \ -\log |z\wedge w| + \Hhat_{F}(z) + \Hhat_{F}(w)
- r(F),
\]
where $r(F) = \frac{1}{d(d-1)}\log |\Res(F)|$.

An explanation for the appearance of the constant term $-r(F)$ will be given shortly.

\medskip

According to \cite{BR}, the function
$G_{F}$ is doubly scale-invariant, in the sense
that if $\alpha, \beta \in K^*$, then
\[
G_{F}(\alpha z, \beta w) \ = \ G_{F}(z,w).
\]

In addition, for all $\gamma \in K^*$, we have
\[
G_{\gamma F}(z, w) \ = \ G_{F}(z,w).
\]

In particular, $G_{F}$ descends to a well-defined
function $g_{\varphi}(z,w)$ on $\PP^1(K)$:
for $z, w \in \PP^1(K)$ and any lifts $\tz, \tw \in K^2$,
\begin{equation} 
\label{gFvDef} 
g_{\varphi}(z,w) \ = \ -\log |\tz\wedge \tw| + \Hhat_{F}(\tz) 
                      + \Hhat_{F}(\tw) - r(F).
\end{equation}
If $z\neq w$ then the right-hand side of (\ref{gFvDef}) is finite; if
$z=w$ then we define $g_{\varphi}(z,z) := +\infty$.

\begin{remark}
In \cite{BR}, it is proved in the case $K=\CC_v$ that the function $g_{\varphi}(z,w)$ is a normalized Arakelov-Green's function
for the canonical measure $\mu_{\varphi}$ on the Berkovich analytic space $\PP^1_{\Berk}$ over $\CC_v$; we refer to
\cite{BR} for a precise explanation of these objects.  Briefly, the fact that $g_{\varphi}(z,w)$ is a normalized Arakelov-Green's function
means that it is the unique solution to the differential equation
\begin{equation}
\label{DiffEQ}
\Delta_z g_{\varphi}(z,w) = \delta_w - \mu_{\varphi}
\end{equation}
subject to the normalization condition
\begin{equation}
\label{Normalization}
\iint g_{\varphi}(z,w) d\mu_{\varphi}(z) d\mu_{\varphi}(w) = 0.
\end{equation}
When $v$ is archimedean, $\PP^1_{\Berk} = \PP^1(\CC)$ and $\mu_{\varphi}$ is the canonical probability measure on 
$\PP^1(\CC)$ introduced by Lyubich and Freire-Lopes-Ma{\~n}{\'e} whose support is equal to the Julia set of $\varphi$.
For our purposes, one can take (\ref{DiffEQ}) as the definition of $\mu_{\varphi}$.
The constant $r(F) = \frac{1}{d(d-1)}\log |\Res(F)|$ which appears in the definition (\ref{gFvDef}) of 
$g_{\varphi}$ is chosen precisely so that (\ref{Normalization}) holds; see \cite[Corollary 4.13]{BR}.
\end{remark}

\begin{remark}
\label{GoodReductionRemark}
If $K$ is non-archimedean and $\varphi$ has good reduction, then $g_\varphi$ is nonnegative and 
Theorem~\ref{GreenTheorem} is trivial.
\end{remark}

\subsection{The homogeneous transfinite diameter of the filled Julia set}

In order to explain some of the applications of (\ref{mainboundequation}), 
we recall some additional facts from \cite{BR}.

Let $E \subset K^2$ be a bounded set.
By analogy with the classical transfinite diameter, define
\[
d^0_n(E) \ := \ \sup_{z_1,\ldots,z_n \in E}
\left( \prod_{i \neq j} |z_i \wedge z_j|  \right)^{\frac{1}{n(n-1)}}.
\]

By \cite{BR}, the sequence of nonnegative real numbers $d^0_n(E)$ is
non-increasing.  In particular, the quantity $d^0_\infty(E) :=
\lim_{n\to\infty} d^0_n(E)$ is well-defined.
We call $d^0_\infty(E)$ the {\em homogeneous transfinite diameter} of
$E$.

The {\emph{filled Julia set}} $K_{F}$ of $F$ in $K^2$
is the set of all $z \in K^2$ for which the iterates $F^{(n)}(z)$
remain bounded.  
Clearly $F^{-1}(K_{F}) = K_{F}$, and the same is true for each $F^{(-n)}$.
It is proved in \cite{BR} that $K_{F}$ is a closed and bounded subset of $K^2$, and that
\[
K_{F} = \{ z \in K^2 \; : \; \Hhat_{F}(z) \leq 0 \}.
\]

The following theorem from \cite{BR}
is a generalization from $K=\CC$ to an arbitrary valued field
of a formula of DeMarco \cite{DeMarco}:

\begin{theorem}
\label{AdelicDeMarcoCapacityTheorem}
If $K$ is algebraically closed, then
\begin{equation}
\label{AdelicDeMarcoFormula}
d^0_\infty(K_{F}) \ = \ |\Res(F)|^{-1/d(d-1)}.
\end{equation}
\end{theorem}

The proof of this result given in \cite{BR} is rather indirect, and requires the development of a lot of capacity-theoretic
machinery, including a detailed analysis of the relationship between the homogeneous transfinite diameter and the
homogeneous sectional capacity.  The principal application of Theorem~\ref{AdelicDeMarcoCapacityTheorem} given in \cite{BR} 
is an adelic equidistribution theorem for points of small dynamical height with respect to a rational map defined
over a number field.  This application requires only the inequality
\begin{equation}
\label{TransfiniteDiameterInequality}
d^0_\infty(K_{F}) \ \leq \ |\Res(F)|^{-1/d(d-1)}.
\end{equation}

Theorem~\ref{GreenTheorem}
furnishes a simpler and more direct proof of
(\ref{TransfiniteDiameterInequality}) 
(which does not require that the field $K$ be algebraically closed).
Indeed, if $z_1,\ldots,z_N$ 
are chosen to belong to $K_{F}$, then $\Hhat_{F}(z_i) \leq 0$ and
(\ref{mainboundequation}) yields
\[
\frac{1}{N(N-1)}\sum_{i \neq j} \log |z_i \wedge z_j| \leq C \frac{\log N}{N-1} - r(F).
\]
Passing to the limit as $N$ tends to infinity and exponentiating gives (\ref{TransfiniteDiameterInequality}).


\subsection{Application to canonical heights}

We now consider a global application of (\ref{mainboundequation}).
Let $k$ be a number field, let $\varphi \in k(T)$ be a rational map with coefficients in $k$, 
and let $\hhat_\varphi$ denote the Call-Silverman 
canonical height function attached to $\varphi$; this can be defined for $P=[z_0:z_1] \in \PP^1(k)$ by
\[
\hhat_\varphi(P) = \sum_{v \in M_k} \frac{[k_v : \QQ_v]}{[k:\QQ]} \Hhat_{F,v}(z_0,z_1). \\
\]
Here $\Hhat_{F,v}$ denotes the function $\Hhat_{F}$ 
over the completion $k_v$ (or over $\CC_v$).

The proof of the following result will be given in \S\ref{ApplicationSection}.

\begin{theorem}
\label{GlobalHeightTheorem}
There exist constants $A,B>0$, depending only on $\varphi$ and $k$, 
such that for any finite extension $k'/k$ with $D = [k':\QQ]$, we have 
\begin{equation}
\label{GlobalHeightBound}
\# \left\{ P \in \PP^1(k') \; : \; \hhat_\varphi(P) \leq \frac{A}{D} \right\} \leq BD\log D.
\end{equation}
\end{theorem}

\begin{remark}
In the special case $\varphi(T)=T^2$, the function $\hhat_\varphi$ is the standard logarithmic Weil height,
and (\ref{GlobalHeightBound}) implies that there is a constant $C>0$ such that if $\alpha \in (k')^*$ has $h(\alpha)>0$ (i.e., 
$\alpha$ is not a root of unity), then
\begin{equation}
\label{WeakLehmer}
h(\alpha) \geq \frac{C}{D^2 \log D}.
\end{equation}
Indeed, if $M$ is the largest integer such that $(M-1)h(\alpha) = h(\alpha^{M-1}) \leq \frac{A}{D}$, then
$\{ 1, \alpha, \alpha^2, \ldots, \alpha^{M-1} \} \subseteq \{ P \in \PP^1(k') \; : \; h(P) \leq \frac{A}{D} \}$,
and thus $M \leq B \cdot D\log D$.  By the maximality of $M$, we have $Mh(\alpha) > \frac{A}{D}$, and therefore
$h(\alpha) > \frac{A}{B} \cdot \frac{1}{D^2 \log D}$.
Of course, much stronger results are known in connection with Lehmer's problem (see e.g. \cite{Dobrowolski}).
Similarly, applying (\ref{GlobalHeightBound}) to a Latt{\`e}s map (for which $\hhat_\varphi$ is the N{\'e}ron-Tate canonical height on the
elliptic curve $E/k$) and using the quadraticity of the canonical height gives the bound
\begin{equation}
\label{WeakEllipticLehmer}
\hhat_E(P) \geq \frac{C}{D^3 \log^2 D}
\end{equation}
for all non-torsion points $P \in E(k')$.
This is precisely the same bound for the elliptic Lehmer problem obtained by Masser \cite{Masser}.
Unfortunately, for a general rational map, Theorem~\ref{GlobalHeightTheorem} does not seem to imply
a Lehmer-type estimate such as (\ref{WeakLehmer}) or (\ref{WeakEllipticLehmer}) which is polynomial in $1/D$.
\end{remark}

Finally, we note that one can deduce from Theorem~\ref{GlobalHeightTheorem} the following result concerning 
fields of definition for preperiodic points of $\varphi$:

\begin{cor}
\label{DegreeCor}
There exists a constant $C$ such that if $P_1,\ldots,P_N$ are distinct preperiodic points of $\varphi$ defined over $k'$, 
then 
\[
[k':k] \geq C \frac{N}{\log N}.
\]
In particular, if $k_n$ denotes the extension of $k$ obtained by adjoining to $k$ the $N_n$ points of exact period $n$,
then
\begin{equation}
\label{DegreeInequality}
[k_n:k] \gg \frac{N_n}{\log N_n}.
\end{equation}
\end{cor}

\begin{remark}
For a generic rational map $\varphi$, we would expect that something stronger than
Corollary~\ref{DegreeCor} should hold.  On the other hand, there are rational maps for
which Corollary~\ref{DegreeCor} is nearly sharp.  For example, if $\varphi(T)=T^2$ and $k=\QQ$,
then the $N$th roots of unity are distinct preperiodic points of $\varphi$ defined over $\QQ(\zeta_N)$
and $[\QQ(\zeta_N):\QQ] = \phi(N)$.  Since $N/\log\log N \ll \phi(N) \ll N$, Corollary~\ref{DegreeCor} is
not far from the truth in this case.
\end{remark}



\section{Discussion and proof of the main result}

\subsection{A sharpening of Theorem~\ref{GreenTheorem}}
\label{SharpeningSection}

Let $(x_1,y_1),\ldots,(x_N,y_N)$ be
nonzero points in the filled Julia set $K_{F}$
whose images in $\PP^1(K)$ are distinct.
Our proof of (\ref{mainboundequation}) comes from a lower bound for the sum
\[
\sum_{i \neq j} -\log |(x_i,y_i) \wedge (x_j,y_j)|
\]
when $N$ belongs to the subset 
\[
\Sigma = \{ N \in \NN \; | \; N = td^k \textrm{\; for some integers \;} 2\leq t \leq 2d-1 \textrm{\; and \; } k\geq 1\}
\]
of the natural numbers.

Specifically, we will prove the following technical result.  For notational convenience, let 
$\epsilon_K$ be zero if the absolute value on $K$ is non-archimedean, and $1$ if it is archimedean.

\begin{theorem}
\label{MainTechnicalTheorem}
Let $N = td^k \in \Sigma$, and let $z_1,\ldots,z_N$ be nonzero 
elements of the filled Julia set $K_{F}$
whose images in $\PP^1(K)$ are all distinct.
Let $R(F)$ denote the smallest radius so that 
$K_{F} \subseteq \{ z \in K^2 \; : \; \| z \| \leq R(F) \}$.
Then
\[
\begin{array}{lll}
\sum_{i \neq j} -\log |z_i \wedge z_j| 
& \geq & r(F) N^2  -\epsilon_K N\log N \\
& & - 2\left( \log R(F) \right) \alpha N  - r(F)(1+\alpha) N, \\
\end{array}
\]
where 
$\alpha = t-1 + (d-1)k > 0$ satisfies
$2 \leq \alpha \leq (d-1)(\log_d N + 2)$.

\end{theorem}

\begin{remark}
Using the fact that $|T z_i \wedge Tz_j| = |z_i \wedge z_j|$ for all $T \in \GL_2(K)$ with
$\det(T)=\pm 1$, it can be shown that the quantity $R(F)$ appearing
in the statement of Theorem~\ref{MainTechnicalTheorem}
may be replaced by the smallest radius $R$ such
that $K_{F} \subseteq \{ z \; : \; \| Tz \| \leq R \}$ for some $T$ with
$\det(T)=\pm 1$.
\end{remark}

As a corollary of Theorem~\ref{MainTechnicalTheorem}, we obtain:

\begin{cor}
\label{GreensFunctionCorollary}
Keeping the notation of Theorem~\ref{MainTechnicalTheorem}, let $N \in \Sigma$,
let $z_1,\ldots,z_N$ be distinct points in $\PP^1(K)$, and set
\[
D_\varphi(z_1,\ldots,z_N) := \sum_{i \neq j} g_{\varphi}(z_i,z_j).
\]
Then
\begin{equation}
\label{GFCeq}
D_\varphi(z_1,\ldots,z_N) \geq  -\epsilon_K N\log N - \left( 2\log R(F)+r(F) \right) \alpha N.
\end{equation}
In particular, 
\[
\liminf_{N \to \infty} \inf_{z_1,\ldots,z_N \in \PP^1(K)} \frac{1}{N(N-1)} 
\sum_{i \neq j} g_{\varphi}(z_i,z_j) \geq 0.
\]
\end{cor}

\begin{remark}
If $K$ is non-archimedean and $\varphi$ has good reduction, then
$\epsilon_K = \log R(F) = r(F) = 0$, and (\ref{GFCeq}) just says that
$D_\varphi(z_1,\ldots,z_N) \geq 0$, which of course already follows from
Remark~\ref{GoodReductionRemark}.
\end{remark}

\subsection{Proof of Theorem~\ref{MainTechnicalTheorem}}

An outline of the proof of Theorem~\ref{MainTechnicalTheorem} is as follows.  
First, we express $\prod_{i \neq j} |(x_i,y_i) \wedge (x_j,y_j)|$ as the determinant of a 
Vandermonde matrix $S$.  We then replace this matrix with 
a new matrix $H$ whose entries involve $F_1^{(k)}(x_i,y_i)$ and $F_2^{(k)}(x_i,y_i)$ for various $k \geq 0$ rather than just 
the standard monomials $x_i^a y_i^b$. 
Replacing $S$ by $H$ amounts to choosing a different basis for the space of homogeneous polynomials in $x$ and $y$ of degree $N-1$, and
to calculate the determinant of the change of basis matrix 
we use a generalization of Lemma 6.5 of \cite{BR}.  
This is the key step in the argument, and is the place where we use the hypotheses that $N \in \Sigma$.  
Finally, we use Hadamard's inequality to estimate the determinant of $H$, using the fact that
$\| F^{(k)}(x_i,y_i) \| \leq R(F)$ for all $k \geq 0$.

\medskip

Let $\Gamma^0(m)$ denote the vector space of homogeneous polynomials of degree $m$ in $K[x,y]$, which has dimension
$N = m+1$ over $K$.  If $N \in \Sigma$, i.e., if $m = td^k - 1$ with $2 \leq t\leq 2d-1$ and $k \geq 1$, we consider the
collection $H(m)$ of polynomials
\[
\begin{array}{lll}
H(m) & = & \{ x^{a_0}y^{b_0} F_1(x,y)^{a_1} F_2(x,y)^{b_1} \cdots F_1^{(k)}(x,y)^{a_k} F_2^{(k)}(x,y)^{b_k} \; | \; \\
& & a_i + b_i = d-1 \textrm{\; for \; } 0\leq i\leq k-1 \textrm{ \; and \; } a_k + b_k = t-1 \} \\ 
& \subset & \Gamma^0(m). \\
\end{array}
\]

The cardinality of $H(m)$ is easily seen to be $N = \dim \Gamma^0(m)$.  
The following proposition shows that $H(m)$ forms a basis for
$\Gamma^0(m)$, and explicitly calculates the determinant of the change of basis matrix between $H(m)$ and
the standard monomial basis $S(m)$ given by
\[
S(m) = \{ x^a y^b \; | \; a+b = m \}.
\]
Its proof is modeled after Lemma 6.5 of \cite{BR}.  

\begin{proposition}
\label{ChangeOfBasisProp}
Let $A$ be the matrix expressing the polynomials $H(m)$ (in some order) in terms of some ordering of the standard basis $S(m)$.
Then $\det(A) = \pm \Res(F)^r$, where
\[
r = \frac{N^2}{2d(d-1)} - \frac{N}{2d(d-1)} \left( t + k(d-1) \right).
\]
In particular, since $\Res(F) \neq 0$, $H(m)$ is a basis for $\Gamma^0(m)$.
\end{proposition}

\begin{proof}
Let $H_1,\ldots,H_N$ be an ordering of the elements of $H(m)$, and let $S_1,\ldots,S_N$ be an ordering of
the elements of $S(m)$, so that $A$ is the $N \times N$ matrix whose $(i,j)$th entry is the coefficient of the monomial $S_i$ 
which appears in the expansion of $H_j$ as a polynomial in $x$ and $y$.
We have $\det(A) = 0$ if and only if some nontrivial linear combination of the elements of $H(m)$ is zero.

Suppose $\det(A) = 0$.  
Then there exist homogeneous polynomials $h_1 \in H(d^k - 1)$ and 
$h_2 \in H((t-1)d^k - 1)$, not both zero, such that 
\[
h_1 \left( F_1^{(k)}(x,y) \right)^{t-1} + h_2 F_2^{(k)}(x,y) = 0.
\]
We may assume that neither of $h_1,h_2$ is the zero polynomial.  
Thus $F_2^{(k)}$ divides $h_1 \left( F_1^{(k)} \right)^{t-1}$.
Since $\deg(h_1) < \deg(F_2^{(k)})$, it follows that
$F_1^{(k)}$ and $F_2^{(k)}$ have a common irreducible factor.  
Thus $\Res(F^{(k)}) = 0$.  But $\Res(F^{(k)})$ is a power of $\Res(F)$ by \cite[Corollary 6.4]{DeMarco},
so $\Res(F) = 0$ as well.  

Conversely, suppose $\Res(F)=0$.  
Then by a standard fact about resultants (see \cite[\S 5.8]{VdW}),
there is a nontrivial relation of the form
\begin{equation}
\label{NontrivialRelation}
h_1 F_1 + h_2 F_2 = 0
\end{equation}
with $h_1,h_2 \in K[x,y]$ homogeneous of degree $d-1$.
If $k = 1$, this already implies that $\det(A)=0$.  If $k\geq 2$, then 
multiplying both sides of (\ref{NontrivialRelation}) by $G(x,y)$, where
\[
G(x,y) = F_1^{d-2} (F_1^{(2)})^{d-1} \cdots (F_1^{(k)})^{t-1},
\]
gives a linear relation which shows that $\det(A)=0$.

\medskip

Now expand both $\det(A)$ and $\Res(F)$ as polynomials in the coefficients of $F_1$ and $F_2$.
Since $\Res(F)$ is irreducible (see \cite[\S 5.9]{VdW}), we find that
\[
\det(A) = C \cdot \Res(F)^r
\]
for some $C \in K^*$ and some natural number $r$.
Now $\Res(F)$ is homogeneous of degree $2d$ in the coefficients of $F_1$ and $F_2$,
and a straightforward calculation shows that $F_1^{(j)}$ and $F_2^{(j)}$ are each homogeneous
of degree $(d^j - 1)/(d-1)$ in the coefficients of $F_1$ and $F_2$.
It follows that $\det(A)$ is homogeneous of degree 
\[
r' = N \left( \sum_{j=1}^{k-1} (d^j - 1) + (t-1)\frac{d^k - 1}{d-1} \right)
\]
in the coefficients of $F_1$ and $F_2$.
Comparing degrees and performing some straightforward algebraic manipulations, we find that 
\[
r = \frac{r'}{2d} = \frac{N^2}{2d(d-1)} - \frac{N}{2d} \left( \frac{t}{d-1} + k \right).
\]

Finally, to compute $C$ we set $F_1 = x^d$ and $F_2 = y^d$, in which case
$H(m)$ is just a permutation of the standard monomial basis $S(m)$.
It follows that $C = \pm 1$.
\end{proof}

Before turning to the proof of Theorem~\ref{MainTechnicalTheorem}, we first recall Hadamard's
inequality.  If $v = (v_1,\ldots,v_N)^{\rm T} \in K^N$, define $\| v \|$ to be 
the $L^2$-norm $\| v \| = (\sum |v_i|^2)^{1/2}$ if $K$ is archimedean, and to be the
sup-norm $\| v \| = \sup |v_i|$ if $K$ is non-archimedean.  If $H$ is a matrix with columns
$h_1,\ldots,h_N \in K^N$, then Hadamard's inequality states that
\begin{equation}
\label{HadamardInequality}
|\det(H)| \leq \prod_{i=1}^N \| h_i \|.
\end{equation}
If $K$ is non-archimedean, (\ref{HadamardInequality}) is immediate from the definition of the determinant and the
ultrametric inequality.  If $K = \RR$ or $\CC$ with the usual absolute value, (\ref{HadamardInequality}) is a well-known
result from linear algebra.  Finally, (\ref{HadamardInequality}) for general archimedean $K$ can be deduced from the cases
$K=\RR,\CC$ using a theorem of Ostrowski (see \cite[Chapter II, Theorem 4.2]{Neukirch}) which says that a complete archimedean
valued field $K$ is isometric to either $(\RR,|\, |^s)$ or $(\CC,| \, |^s)$ for some $0 < s  \leq 1$.

\medskip

We can now give the proof of Theorem~\ref{MainTechnicalTheorem}.

\begin{proof}
Let $\{ S_1,\ldots, S_N \}$ and $\{ H_1,\ldots,H_N \}$ be as in the proof of 
Proposition~\ref{ChangeOfBasisProp}.
Using the homogeneous version of the standard formula for the determinant of a Vandermonde matrix, we have
\begin{equation}
\label{VandermondeIdentity}
\prod_{i \neq j} |x_iy_j - x_j y_i| = |\det(S)|^2,
\end{equation}
where $S$ is the matrix whose $(i,j)$th entry is $S_j(x_i,y_i)$.
If $H$ is the matrix whose $(i,j)$th entry is $H_j(x_i,y_i)$, then $H = SA$,
so that by Proposition~\ref{ChangeOfBasisProp}, we have 
\begin{equation}
\label{ResultantFormula}
\det(S)^2 = \det(H)^2 (\det(A))^{-2} = \det(H)^2 \Res(F)^{-2r}.
\end{equation}

On the other hand, we can estimate $|\det(H)|^2$ using 
(\ref{HadamardInequality}).
Letting $h_i$ be the $i$th column of $H$, we obtain
\begin{equation}
\label{HadamardBound}
\begin{array}{lll}
|\det(H)|^2 & \leq & \prod_{i=1}^N \| h_i \| \\
& \leq & N^{\epsilon_K N} \prod_i R(F)^{2\left( k(d-1)+(t-1)\right)} \\
& = & N^{\epsilon_K N} R(F)^{\left( 2t-2 + k(2d-2)\right) N}. \\
\end{array}
\end{equation}

Putting together (\ref{VandermondeIdentity}), (\ref{ResultantFormula}), and
(\ref{HadamardBound}) gives
\begin{equation}
\label{UpperBound}
\begin{array}{lll}
\prod_{i \neq j} |x_iy_j - x_j y_i| & = & |\det(S)|^2 \\
& = & |\det(H)|^2 \cdot |\Res(F)|^{-2r} \\
& \leq & N^{\epsilon_K N} R(F)^{2\alpha N} |\Res(F)|^{-2r}, \\
\end{array}
\end{equation}
where $-2r = - \frac{N^2}{d(d-1)}  + \frac{N}{d(d-1)} \left( t + k(d-1) \right)$.

Taking the negative logarithm of both sides of (\ref{UpperBound}) 
gives the desired lower bound
\[
\begin{array}{lll}
\sum_{i \neq j} -\log |x_iy_j - x_j y_i| & \geq & -\epsilon_K N\log N - 2N\alpha \log R(F) \\
& & + \frac{N^2}{d(d-1)} \log |\Res(F)| \\ 
& & - \frac{N}{d(d-1)} \left( t + k(d-1) \right) \log|\Res(F)|. \\
\end{array}
\]

\end{proof}


Corollary~\ref{GreensFunctionCorollary} is deduced from Theorem~\ref{MainTechnicalTheorem} as follows.

\begin{proof}
Without loss of generality, we may assume that the valuation on $K$ is nontrivial, and we may replace $K$ by $\Kbar$, 
We may therefore assume that the value group $|K^*|$ 
is dense in the group $\RR_{>0}$ of nonnegative reals.
Since $\Hhat_{F}(cz) = \Hhat_{F}(z) + \log |c|$ for all $z \in K^2 - \{ 0 \}$ and all $c \in K^*$,
given $\varepsilon > 0$ we can choose coordinates $(x_i,y_i)$ for $z_i$ so that $-\varepsilon \leq \Hhat_{F}(x_i,y_i) \leq 0$.
Then $(x_i,y_i) \in K_{F}$ and we can apply Theorem~\ref{MainTechnicalTheorem} to $(x_1,y_1),\ldots,(x_N,y_N)$.
Simplifying the resulting expression and letting $\varepsilon \to 0$ gives the desired inequality.
\end{proof}

Finally, we explain how to deduce Theorem~\ref{GreenTheorem} from
Corollary~\ref{GreensFunctionCorollary}.
We first recall the following lemma from \cite[Lemma 3.48]{BR}:

\begin{lemma}
\label{monotoniclemma}
For $N\geq 2$, define
\[
D_N := \inf_{z_1,\ldots,z_N \in \PP^1(K)} \frac{1}{N(N-1)} D_\varphi(z_1,\ldots,z_N).
\]
Then the sequence $D_N$ is non-decreasing.
\end{lemma}

We now give the proof of Theorem~\ref{GreenTheorem}:

\begin{proof}

Without loss of generality, we may assume that $N\geq 2d$.  Let $N'$ be the smallest integer less than
or equal to $N$ which belongs to $\Sigma$.  One deduces easily from the definition of $\Sigma$ that
\begin{equation}
\label{NearestGoodNIneq}
\frac{N-1}{2} \leq N' - 1 \leq N - 1.
\end{equation}
Since $N' \in \Sigma$, it follows from Corollary~\ref{GreensFunctionCorollary} that there is a constant 
$C'>0$ (independent of $N$ and $N'$) such that 
\[
D_{N'} \geq -C' \frac{\log N'}{N'-1}.
\]
From this, Lemma~\ref{monotoniclemma}, and (\ref{NearestGoodNIneq}), it follows that
\[
\frac{1}{N(N-1)} D_\varphi(z_1,\ldots,z_N) \geq D_{N'} \geq 
-C \frac{\log N}{N-1}
\]
where $C = 2C'$.
This immediately yields the desired result.
\end{proof}


\section{Global applications}
\label{ApplicationSection}

In this section, we give a proof of Theorem~\ref{GlobalHeightTheorem} and Corollary~\ref{DegreeCor}.
The idea is to use a pigeonhole principle argument as in \cite{HindrySilvermanI},\cite{HindrySilvermanII},
\cite{HindrySilvermanIII}.  
Throughout this section, $k$ denotes a number field and $\varphi \in k(T)$ is a rational function of
degree at least 2.  We will denote by $M_k$ the set of places of $k$, and by $g_{\varphi,v}$ the 
normalized Arakelov-Green's function over $k_v$ (or $\CC_v$) attached to $\varphi$.
We begin with the following lemma:

\begin{lemma}
\label{positivitylemma}
Suppose $v \in M_k$ is archimedean.  Then:
\begin{itemize}
\item[(a)] The function
\[
g_{\varphi,v}(z,w) : \PP^1(\CC) \times \PP^1(\CC) \to \RR \cup \{ \infty \}
\]
is lower semicontinuous.
\item[(b)]
Given $\delta > 0$,
there exists a finite covering $U_1,\ldots,U_s$ of $\PP^1(\CC)$ by open subsets so that for each 
$1 \leq i \leq s$, $g_{\varphi,v}(z,w) > \delta$ for all $z,w \in U_i$.
\end{itemize}
\end{lemma}

\begin{proof}
(a) It suffices to note that $g_{\varphi,v}$ is an increasing limit of the continuous functions
\[
g^{(M)}_{\varphi,v} = \max(g_{\varphi,v},M)
\]
for $M > 0$ (see \cite[\S 3.5]{BR}).

(b) By part (a), the fact that $g_{\varphi,v}=+\infty$ along the diagonal,
and the definition of the product topology on $\PP^1(\CC) \times \PP^1(\CC)$,
for each $x \in \PP^1(\CC)$ there is an open neighborhood
$U_x$ of $x$ in $\PP^1(\CC)$ such that $g_{\varphi,v}(z,w) > \delta$ for $z,w \in U_x$.
By compactness of $\PP^1(\CC)$, the covering $\{ U_x \; | \; x \in \PP^1(\CC) \}$ has a finite
subcover $U_1,\ldots,U_s$.
\end{proof}

We can now prove Theorem~\ref{GlobalHeightTheorem}.

\begin{proof}
Fix a finite extension $k'/k$ and let $D = [k':\QQ]$.  
By the product formula, for all $z,w \in \PP^1(k')$ with $z\neq w$ we have
\begin{equation}
\label{globaleq1}
\sum_{v \in M_{k'}} 
\frac{[k'_v : \QQ_v]}{[k':\QQ]} g_{\varphi,v}(z,w) = \hhat_{\varphi}(z) + \hhat_{\varphi}(w).
\end{equation}
Let $v_0$ be a fixed archimedean place of $k'$; then by Lemma~\ref{positivitylemma}, there exists
an open covering $U_1,\ldots,U_s$ of $\PP^1(\CC)$ and a constant $C_1 > 0$ such that
$g_{\varphi,v_0}(z,w) \geq C_1$ whenever $z,w \in U_i$ ($i=1,\ldots,s$).
If $z_1,\ldots,z_N \in \PP^1(k')$, then 
by Theorem~\ref{GreenTheorem} and the fact that $g_{\varphi,v} \geq 0$ whenever 
$\varphi$ has good reduction at $v$,
there is a constant $C_2>0$ depending only on $\varphi$ and $k$ such that
\[
g(z_1,\ldots,z_N) := \sum_{v \in M_{k'}} \sum_{i \neq j} \frac{[k'_v : \QQ_v]}{[k':\QQ]} g_v(z_i,z_j) \geq -C_2 N\log N.
\]
Moreover, if $z_1,\ldots,z_N \in U_i$ for some $i$, then
\begin{equation}
\label{globaleq2}
g(z_1\ldots,z_N) \geq \frac{C_1}{D} N^2 -C_2 N\log N.
\end{equation}
Let $M = [\frac{N-1}{s}]+1$.  By the pigeonhole principle, in any subset of $\PP^1(\CC)$ of cardinality $N$, there is an $M$-element subset contained in
some $U_i$.  Without loss of generality, order the $z_j$'s so that this subset is $\{ z_1,\ldots,z_{M} \}$.  
Applying (\ref{globaleq1}) and (\ref{globaleq2}), we obtain
\[
\frac{C_1}{D} M^2 - C_2 M\log M \leq g(z_1,\ldots,z_{M}) \leq 2M^2 \max_j \hhat_{\varphi}(z_j).
\]
If $\hhat_{\varphi}(z_j) \leq \frac{C_1}{4D}$ for all $j=1,\ldots,N$, then we obtain
\[
\frac{C_1}{2D} M \leq C_2 \log M,
\]
which implies that $M \leq B' D \log D$ for some constant $B' > 0$ depending 
only on $\varphi$ and $k$.
Thus $N \leq Ms+1 \leq B D \log D$ for some $B>0$ depending only on $\varphi$ and $k$.
Setting $A = C_1/4$ now gives the desired result.  
\end{proof}

Corollary~\ref{DegreeCor} follows easily:

\begin{proof}
Let $\{ P_1,\ldots,P_N \}$ be a set of preperiodic points defined over a number field $k'$ with
$[k':k] = D$.  Then by Theorem~\ref{GlobalHeightTheorem}, we have
$N \leq B D \log D$, which implies that $D \gg N/\log N$ as claimed.
\end{proof}

\end{document}